\documentclass[a4paper,11pt]{article}
\usepackage[latin1]{inputenc}
\usepackage[francais,english]{babel}    
\usepackage{amssymb}
\usepackage{amsmath}
\usepackage{amsthm}
\usepackage{multicol}
\usepackage{graphicx}
\usepackage{color}
\usepackage[T1]{fontenc}
\usepackage{yfonts}
\usepackage{placeins}
\usepackage{array}
\usepackage{dsfont}
\usepackage{stmaryrd}
\usepackage{verbatim}
\usepackage{caption}
\usepackage{transparent}
\usepackage{bbold}
\usepackage{hyperref}
\usepackage{enumitem}
\usepackage{tikz}
\usetikzlibrary{patterns}
\usepackage{tikz-3dplot}

\newtheorem{thm}{Theorem}[section]
\newcommand{\bthm}{\begin{thm}}
\newcommand{\ethm}{\end{thm}}

\newtheorem{thmi}{Theorem}
\newcommand{\bthmi}{\begin{thmi}}
\newcommand{\ethmi}{\end{thmi}}

\newtheorem{cori}[thmi]{Corollary}
\newcommand{\bcori}{\begin{cori}}
\newcommand{\ecori}{\end{cori}}

\newtheorem{mthm}{Theorem}
\newcommand{\bmthm}{\begin{mthm}}
\newcommand{\emthm}{\end{mthm}}

\newtheorem{mcor}[mthm]{Corollary}
\newcommand{\bmcor}{\begin{mcor}}
\newcommand{\emcor}{\end{mcor}}
\newtheorem{mconj}[mthm]{Conjecture}
\newcommand{\bmconj}{\begin{mconj}}
\newcommand{\emconj}{\end{mconj}}

\newtheorem*{conj}{Conjecture}
\newcommand{\bconj}{\begin{conj}}
\newcommand{\econj}{\end{conj}}

\newtheorem*{question}{Question}
\newcommand{\bq}{\begin{question}}
\newcommand{\eq}{\end{question}}

\newtheorem*{thn}{Theorem}
\newcommand{\bthn}{\begin{thn}}
\newcommand{\ethn}{\end{thn}}

\newtheorem{exo}{Exercise}
\newcommand{\bex}{\begin{exo}}
\newcommand{\eex}{\end{exo}}

\newtheorem{sol}{Solution}
\newcommand{\bsol}{\begin{sol}}
\newcommand{\esol}{\end{sol}}

\newtheorem{pro}[thm]{Proposition}
\newcommand{\bpro}{\begin{pro}}
\newcommand{\epro}{\end{pro}}

\newtheorem{cor}[thm]{Corollary}
\newcommand{\bcor}{\begin{cor}}
\newcommand{\ecor}{\end{cor}}

\newtheorem{lem}[thm]{Lemma}
\newcommand{\blem}{\begin{lem}}
\newcommand{\elem}{\end{lem}}

\theoremstyle{definition}

\newtheorem{defi}[thm]{Definition}
\newcommand{\bdf}{\begin{defi}}
\newcommand{\edf}{\end{defi}}

\newtheorem*{defis}{Definition}
\newcommand{\bdfs}{\begin{defis}}
\newcommand{\edfs}{\end{defis}}

\newtheorem*{rmk}{Remark}
\newcommand{\brk}{\begin{rmk} \upshape}
\newcommand{\erk}{\end{rmk}}

\newtheorem*{rmks}{Remarks}
\newcommand{\brks}{\begin{rmks} \upshape}
\newcommand{\erks}{\end{rmks}}

\newtheorem*{exe}{Example}
\newcommand{\bexe}{\begin{exe} \upshape}
\newcommand{\eexe}{\end{exe}}

\newtheorem*{exes}{Examples}
\newcommand{\bexes}{\begin{exes} \upshape}
\newcommand{\eexes}{\end{exes}}

\newtheorem*{pre}{Proof}
\newcommand{\bp}{\begin{pre} \upshape}
\newcommand{\ep}{\hfill \qed \end{pre}}
\newcommand{\epp}{\end{pre}}

\newcommand{\beq}{\begin{eqnarray*}}
\newcommand{\eeq}{\end{eqnarray*}}

\newcommand{\beqn}{\begin{equation}}
\newcommand{\eeqn}{\end{equation}}

\newcommand{\ben}{\begin{enumerate}}
\newcommand{\een}{\end{enumerate}}
\newcommand{\bit}{\begin{itemize} \renewcommand{\labelitemi}{$\bullet$} \renewcommand{\labelitemii}{$\star$}}
\newcommand{\eit}{\end{itemize}}

\newcommand{\bfg}{
\begin{figure}[H]
\begin{center}}
\newcommand{\efg}{
\end{center}
\end{figure}
\FloatBarrier}

\newcolumntype{M}[1]{>{\raggedright}m{#1}}

\newcommand{\R}{\mathbb{R}}

\newcommand{\Z}{\mathbb{Z}}

\renewcommand{\SS}{\mathbb{S}}

\newcommand{\bs}{\symbol{92}}

\newcommand{\ov}{\overline}

\newcommand{\st}{\, | \,}

\newcommand{\ra}{\rightarrow}

\renewcommand{\geq}{\geqslant}
\renewcommand{\leq}{\leqslant}

\renewcommand{\>}{\rangle}

\newcommand{\mk}{\medskip}

\newcommand{\sign}{\begin{flushright}
Thomas Haettel \\
IMAG, Univ Montpellier, CNRS, France \\
thomas.haettel@umontpellier.fr
\end{flushright}}

\makeatletter
\def\Ddots{\mathinner{\mkern1mu\raise\p@
\vbox{\kern7\p@\hbox{.}}\mkern2mu
\raise4\p@\hbox{.}\mkern2mu\raise7\p@\hbox{.}\mkern1mu}}
\makeatother

\makeatletter
\def\maketitles{%
  \null
  \thispagestyle{empty}%
  \vfill
  \begin{center}\leavevmode
    \normalfont
    {\LARGE \@title\par}%
    \vskip 1.2cm
    {\large \@author\par}%
    \vskip 1.2cm
    {\large \@subtitle\par}%
    \vskip 0.8cm
    {\large \@date\par}%
  \end{center}%
  \vfill
  \null
  \cleardoublepage
  }

\def\date#1{\def\@date{#1}}
\def\author#1{\def\@author{#1}}
\def\title#1{\def\@title{#1}}
\def\subtitle#1{\def\@subtitle{#1}}
\makeatother

\usetikzlibrary{arrows}

\title{Cubulation of some triangle-free Artin groups}
\author{Thomas Haettel}
\date{\today}

\begin{document}

\selectlanguage{english}

\maketitle

\begin{center}
\begin{minipage}{0.8\textwidth}
\textsc{Abstract.} We prove that some classes of triangle-free Artin groups act properly on locally finite, finite-dimensional CAT(0) cube complexes. In particular, this provides the first examples of Artin groups that are properly cubulated but cannot be cocompactly cubulated, even virtually. The existence of such a proper action has many interesting consequences for the group, notably the Haagerup property, and the Baum-Connes conjecture with coefficients.
\end{minipage}
\end{center}

\let\thefootnote\relax

\section*{Introduction}

Artin groups are natural combinatorial generalizations of Artin's braid groups. For every finite simple graph $\Gamma$ with vertex set $S$ and with edges labeled by some integer in $\{2,3,\dots\}$, one associates the Artin group $A(\Gamma)$ with the following presentation:
$$A(\Gamma) = \<S \st \forall \{s,t\} \in \Gamma^{(1)}, w_m(s,t)=w_m(t,s) \mbox{ if the edge $\{s,t\}$ is labeled $m$}\>,$$
where $w_m(s,t)$ is the word $stst \dots$ of length $m$. Note that when $m=2$, then $s$ and $t$ commute, and when $m=3$, then $s$ and $t$ satisfy the classical braid relation $sts=tst$. Also note that when adding the relation $s^2=1$ for every $s \in S$, one obtains the Coxeter group $W(\Gamma)$ associated to $\Gamma$.

\mk

Apart from certain very particular classes, Artin groups remainly largely mysterious (see~\cite{charney_problems} and \cite{mccammond_mysterious}). An approach to understanding Artin groups is to find actions by isometries on metric spaces with some nonpositive curvature, which usually have very nice algebraic consequences on the groups.

\mk

In order to state what is known, let us recall the following classes of Artin groups. An Artin group $A(\Gamma)$ is called:
\bit
\item \emph{right-angled} if all labels are equal to $2$,
\item of \emph{spherical type} if $W(\Gamma)$ is finite,
\item of \emph{type FC} if every complete subgraph of $\Gamma$ spans a spherical type Artin subgroup, or
\item of \emph{large type} if all labels are greater or equal to $3$, and
\item of \emph{XXL type} if all labels are greater or equal to $5$.
\eit

\mk

Concerning CAT(0) spaces, R.~Charney asks whether every Artin group acts properly and cocompactly on a CAT(0) space (see~\cite{charney_problems}). Very few cases are known, essentially right-angled Artin groups (see~\cite{charney_davis_salvetti}), groups with few generators (see~\cite{brady_artin_three},\cite{brady_mccammond}, \cite{b6}) and groups with sufficiently large labels (see~\cite{brady_noel_crisp}, \cite{brady_mccammond_artinthree}, \cite{bell_artin}, \cite{haettel_xxl}).

\mk

Concerning variations on the notion of nonpositive curvature, Bestvina defined a geometric action of Artin groups of spherical type on a simplicial complex with some nonpositive curvature features (see~\cite{bestvina_artin}). More recently, Huang and Osajda proved (see~\cite{huang_osajda}) that every Artin group of almost large type (a class including all Artin groups of large type) acts properly and cocompactly on systolic complexes, which are a combinatorial variation of nonpositive curvature. They also proved (see~\cite{huang_osajda_helly}) that every Artin group of type FC acts geometrically on a Helly graph, which gives rise to classifying spaces with convex geodesic bicombings.

\mk

Concerning CAT(0) cube complexes, we have a conjectural description of all Artin groups acting properly and cocompactly (even virtually) on a CAT(0) cube complex (see~\cite{haettel_artin_cubic}), which shows that such cocompactly cubulated groups look much like right-angled Artin groups.

\mk

However, Haglund and Wise ask whether every Artin group acts properly on a CAT(0) cube complex, but not necessarily cocompactly (see~\cite{haglund_wise}). Apart from the very few examples of cocompactly cubulated Artin groups, no other example was known. The purpose of this article is to give the first examples of properly cubulated Artin groups.

\bmthm \label{thm:main} Assume that $A=A(\Gamma)$ is an Artin group satisfying one of the following.
\begin{enumerate}
\item[(A)] $\Gamma$ has no cycle.
\item[(B)] $\Gamma$ is bipartite and has all labels at least $3$.
\item[(C)] $\Gamma$ has no triangle and has no label $3$.
\een
Then $A$ is the fundamental group of a locally finite, finite-dimensional, non-compact, locally CAT(0) cube complex with finitely many hyperplanes.
\emthm

Note that all these Artin groups are $2$-dimensional, and as such enter the conjectural classification of cocompactly cubulated Artin groups (see~\cite{haettel_artin_cubic}). Most of these are, in fact, not cocompactly cubulated. Here are possibly the simplest interesting examples.

\bmcor \label{maincor:square}
Let $\Gamma$ be a connected, bipartite graph with diameter at least $3$ and with labels at least $3$. Then $A(\Gamma)$ acts properly on a finite-dimensional, locally finite CAT(0) cube complex, but no finite index subgroup of $A(\Gamma)$ acts properly and cocompactly on a CAT(0) cube complex.
\emcor

Among the many consequences of the existence of a proper action on a CAT(0) cube complex, here are the ones that are new for these Artin groups.

\bmcor \label{maincor:csq}
Let $A$ be an Artin group as in Theorem~\ref{thm:main}.
\ben
\item $A$ satisfies the Haagerup property.
\item $A$ is weakly amenable, with Cowling-Haagerup constant $1$.
\item $A$ satisfies the Baum-Connes conjecture with coefficients.
\item $A$ satisfies the Rapid Decay property RD. 
\item $A$ has finite asymptotic dimension.
\een
\emcor

The Haagerup property for a group $A$ asks for the existence of a metrically proper affine action by isometries on a Hilbert space. Apart from right-angled Artin groups and cocompactly cubulated Artin groups (see~\cite{haettel_artin_cubic}), no example was known. The most notable consequence of the Haagerup property is the Baum-Connes conjecture with coefficients.

The weak amenability of a discrete group $A$ asks for the existence of a sequence of finitely supported functions on $A$, converging pointwise to $1$, which are uniformly bounded (by $1$) in the completely bounded multiplier norm (see~\cite{guentner_higson} for details). As for the Haagerup property, apart from right-angled Artin groups and cocompactly cubulated Artin groups (see~\cite{haettel_artin_cubic}), no example was known.

Concerning the Baum-Connes conjecture, the only examples were essentially braid groups (see~\cite{oyono_oyono_braid} and \cite{schick_extension}), some large type Artin groups (see~\cite{ciobanu_holt_rees}) and XXL type Artin groups (see~\cite{haettel_xxl}).

The property RD aks for a polynomial bound on the norm of a convolution operator. One interesting consequence of property RD is that, together with a proper cocompact action on a CAT(0) space, it implies the Baum-Connes conjecture. The only previous examples were the $4$-strand braid group (see~\cite{barre_pichot_rd}) and many large type Artin groups (see~\cite{ciobanu_holt_rees}).

The asymptotic dimension of a group is coarse notion of dimension, whose finiteness for a finitely generated group implies the Novikov conjecture. The previous examples of Artin groups with finite asymptotic dimension were braid groups (see~\cite{bell_fujiwara}), right-angled Artin groups and cocompactly cubulated Artin groups (see~\cite{haettel_artin_cubic}).

\mk

The idea to construct a locally CAT(0) cube complex which has for fundamental group an Artin group $A(\Gamma)$ as in Theorem~\ref{thm:main} is very simple. We will construct, for each edge $\{a,b\}$ in $\Gamma$, a locally CAT(0) cube complex $M_{a,b}$ whose fundamental group is a dihedral Artin group $\<a,b\>$. We will then glue all these cube complexes along subcomplexes corresponding to the generators of $A(\Gamma)$. However, there are many challenges when applying this strategy.

\mk

The first problem is to build a locally CAT(0) cube complex $M_{a,b}$ for a dihedral Artin group $\<a,b \st w_m(a,b)=w_m(b,a)\>$, with two locally convex subcomplexes $M_a$ and $M_b$ with fundamental groups $\<a\>$ and $\<b\>$ respectively. According to~\cite{haettel_artin_cubic}, this cannot be achieved by a compact CAT(0) cube complex, hence the need to consider only proper actions. In the first two cases (A) and (B), we will use the same complex $M_{a,b}$, and in the last case (C) we will use a slightly different complex.

\mk

The second problem is to ensure that, for each such complex $M_{a,b}$, the subcomplexes $M_a$ and $M_b$ do not intersect too much. This is the first role of the assumptions in Theorem~\ref{thm:main}. In the first two cases (A) and (B), we will consider an action of $\<a,b \st w_m(a,b)=w_m(b,a)\>$ on the product of $\R^m$ and a $m$-regular tree. In this case, $M_a$ and $M_b$ only intersect in one edge. In the last case (C), we will consider an action of $\<a,b \st w_m(a,b)=w_m(b,a)\>$ on the product of $\R^m$ and a tree-like CAT(0) square complex. In this case, $M_a$ and $M_b$ only intersect in a vertex.

In each case (A), (B) and (C), the assumptions ensure that we are able to glue all the complexes $M_{ab}$'s, for each edge $\{a,b\}$ of $\Gamma$, into a cube complex $M$ whose fundamental group is the Artin group $A(\Gamma)$.

\mk

The last question is to decide when the complex $M$ is locally CAT(0). In the case (A), the complex $M$ is a tree-like gluing of locally CAT(0) subspaces $M_{ab}$'s along convex subspaces $M_a$'s, which guarantees that $M$ is locally CAT(0). In the cases (B) and (C), we will apply Gromov's flag link condition: since $\Gamma$ has no triangle, if we consider three pairwise intersecting subcomplexes $M_{ab}$'s, they have a small global intersection, which ensures that links of vertices in $M$ are flag.

\mk

\textbf{Outline of the article:} In Section~\ref{sec:euclidean}, we describe the action of each dihedral Artin group on a Euclidean space that will be used in all three cases (A), (B) and (C). In Section~\ref{sec:AB_tree}, we describe the action of each dihedral Artin group on a tree that will be used in the cases (A) and (B). In Section~\ref{sec:AB_glue}, we explain in the cases (A) and (B)  how to glue the dihedral subcomplexes and obtain a cube complex whose fundamental group is the Artin group, and prove the CAT(0) property. In Section~\ref{sec:C_tree}, we describe the action of each dihedral Artin group on a tree-like CAT(0) square complex that will be used in the case (C). In Section~\ref{sec:C_glue}, we explain in the case (C) how to glue the dihedral subcomplexes and obtain a cube complex whose fundamental group is the Artin group, and prove the CAT(0) property. In the last Section~\ref{sec:corollaries}, we give the proofs of the two corollaries.

\mk

\textbf{Acknowledgments:} The author would like to thank Chris Cashen for discussions and an invitation to the University of Vienna, where part of this work was initiated. The author would also like to thank Anthony Genevois, Piotr Przytycki and Hoel Queffelec for many interesting discussions. We also thank the anonymous referee for useful comments.

\section{The Euclidean action} \label{sec:euclidean}

We start by describing an action of dihedral Artin groups on Euclidean spaces, that will be used in all cases $(A)$, $(B)$ and $(C)$.

\bpro \label{pro:euclidean_action}
For each $m \geq 2$, there exists a cocompact, cubical action of the dihedral Artin group $I_2(m)=\<a,b \st w_m(a,b) = w_m(b,a)\>$ on $\Sigma_{a,b} =\R^m$, with the standard cubical structure, and base vertex $x_0=0$, such that:
\ben
\item The action is given by:
$$a \cdot (y_0,\dots,y_{m-1}) = (y_0+1,y_{-1},y_{-2},\dots,y_{-(m-1)}) $$
$$\mbox{ and } b \cdot (y_0,\dots,y_{m-1}) = (y_2,y_1+1,y_0,y_{2-3},y_{2-4},\dots,y_{2-(m-1)}),$$
with indices in $\Z / m\Z$.
\item The element $a$ acts as a translation on the line $\Sigma_a = \R \times \{0\}^{m-1}$ containing $x_0$.
\item The element $b$ acts as a translation on the line $\Sigma_b = \{0\} \times \R \times \{0\}^{m-2}$ containing $x_0$.
\item The lines $\Sigma_a$ and $\Sigma_b$ intersect in $\{x_0\}$.\een
\epro

\bp
Note that the linear part of the action corresponds to the permutation of the $m$ coordinates given by the action of the Coxeter dihedral group on the vertices of the regular $m$-gon.

\mk

Let us compute the actions of the elements $ab$ and $ba$:
$$ab \cdot (y_0, \dots, y_{m-1}) = (y_2+1,y_3,y_4,\dots,y_{(m-2)+2},y_{(m-1)+2}+1)$$
$$ba \cdot (y_0, \dots, y_{m-1}) = (y_{-2},y_{-1}+1,y_0+1,y_1,\dots,y_{(m-1)-2}).$$

When $m=2p$ is even, we deduce the actions of the elements $w_m(a,b)=(ab)^p$ and $w_m(b,a)=(ba)^p$:
$$w_m(a,b) \cdot (y_0, \dots, y_{m-1}) = w_m(b,a) \cdot (y_0, \dots, y_{m-1}) = (y_0+1,y_1+1,\dots,y_{m-1}+1).$$

When $m=2p+1$ is odd, we deduce the action of $(ab)^p$:
$$(ab)^p \cdot (y_0, \dots, y_{m-1}) = (y_{m-1}+1,y_0,y_1+1,y_2+1,\dots,y_{m-2}+1),$$
so we can compute the actions of $w_m(a,b)=(ab)^pa$ and $w_m(b,a)=b(ab)^p$:
$$w_m(a,b) \cdot (y_0, \dots, y_{m-1}) = w_m(b,a) \cdot (y_0, \dots, y_{m-1}) = (y_1+1,y_0+1,y_{-1}+1,\dots,y_{m-2}+1).$$

In each case, this defines an affine cubical action of $I_2(m)$ on $\R^m$.

To see that the action of $I_2(m)$ is cocompact, notice that the pure Artin subgroup, kernel of the morphism $I_2(m) \ra D_{2 \times m}$ (where $D_{2 \times m}$ denotes the Coxeter dihedral group), acts transitively on the cocompact lattice $(2\Z)^m \subset \R^m$.
\ep

We will need the two following technical results in the sequel.

\blem \label{lem:action_separates}
For every $m \geq 2$, for any $n \in \Z \bs \{0\}$, we have the following:
\ben
\item $\big( (ba)^n \cdot \Sigma_a\big) \cap \Sigma_a = \emptyset$.
\item $\big((ba)^n \cdot \Sigma_a\big) \cap \Sigma_b = \emptyset$.
\item $\big(w_m(a,b)^n \cdot \Sigma_a\big) \cap \Sigma_a = \emptyset$.
\item $\big(w_m(a,b)^n \cdot \Sigma_a\big) \cap \Sigma_b = \emptyset$.
\een
\elem

\bp
Assume first that $m=2$. Then the action of $ab=ba=w_m(a,b)=w_m(b,a)$ on $\Sigma_{a,b}$ is given by $ab \cdot (y_0,y_1)=(y_0+1,y_1+1)$. The result is clear in this case.

\mk

Assume now that $m \geq 3$.

\mk

Recall that the action of $ba$ on $\Sigma_{a,b}$ is given by $ba \cdot (y_0, \dots, y_{m-1}) = (y_{-2},y_{-1}+1,y_0+1,y_1,\dots,y_{(m-1)-2})$. In particular, if $y \in \Sigma_a \cup \Sigma_b$, then for any $n \in \Z \bs \{0\}$, we know that there exists $i \in \{2,3,\dots,m-1\}$ such that $((ba)^n \cdot y)_i \neq 0$, so $(ba)^n \cdot y \not\in \Sigma_a \cup \Sigma_b$.

\mk

Recall that the action of $w_m(a,b)$ on $\Sigma_{a,b}$ has for translation part $(1,1,\dots,1)$. Therefore, if $y \in \Sigma_a \cup \Sigma_b$, then for any $n \in \Z \bs \{0\}$, we know that for every $i \in \{2,3,\dots,m-1\}$ we have $(w_m(a,b)^n \cdot y)_i \neq 0$, so $w_m(a,b)^n \cdot y \not\in \Sigma_a \cup \Sigma_b$.
\ep

\blem \label{lem:action_free}
For every $m \geq 2$, for any $n \in \Z \bs \{0\}$ and for any $y \in \Sigma_{a,b}$, we have the following:
\ben
\item $(ba)^n \cdot y \neq y$.
\item $w_m(a,b)^n \cdot y \neq y$.
\een
\elem

\bp
Notice that if the sum of coordinates of $y$ is $s$, then the sum of coordinates of $(ba) \cdot y$ is $s+2$, and the sum of coordinates of $w_m(a,b) \cdot y$ is $s+m$. The result follows.
\ep

Without risk of confusion, if $a$ is a generator of an infinite cyclic group $\<a\>$, we will denote by $\Sigma_a \simeq \R$ the real line tiled by unit segments, with base vertex $x_0=0$. We also define the action of $\<a\>$ on $\Sigma_a$ where $a$ acts as a translation by $+1$.

\section{The tree actions - cases (A) and (B)} \label{sec:AB_tree}

We will now describe an action by isometries of each dihedral Artin group on a regular tree, that will be used in the cases (A) and (B).

\blem \label{lem:tree_action_A}
For every $m \geq 2$, consider distinct colors $\chi(a),\chi(b)$ in $\{0,1\}$ on the generators of the dihedral Artin group $I_2(m)=\<a,b \st w_m(a,b)=w_m(b,a)\>$. There exists a vertex-transitive action of $I_2(m)$ on an $m$-regular tree $T_{a,b}$ such that the following hold.
\ben
\item The elements $a,b$ act as translations on $T_{a,b}$, with axes $T_a,T_b$ both containing an edge $[t_0,t_1]$.
\item The stabilizer of $t_0$ is $\<ba\>$, and the stabilizer of $t_1$ is $\<ab\>$.
\item If $m \geq 3$, then $T_a \cap T_b = [t_0,t_1]$. If $m=2$, then $T_a=T_b=T_{a,b}$.
\item Let $g \in \{a,b\}$. If $\chi(g)=0$ then $g \cdot t_0=t_1$, and if $\chi(g)=1$ then $g \cdot t_1 = t_0$.
\een
\elem

\bp
Without loss of generality, assume that $\chi(a)=0$ and $\chi(b)=1$.

\mk

Assume first that $m=2$. Let us denote $T_{a,b}=\R$, with the standard tiling by unit segments, with base vertex $t_0=0$ and its neighbour $t_1=1$. Consider the action of $\<a,b\> \simeq \Z^2$ on $T_{a,b}$, where for each $g \in \{a,b\}$, the element $g$ acts on $T_{a,b}$ by a translation of $(-1)^{\chi(g)}$. The axes $T_a,T_b$ of $a,b$ satisfy $T_a=T_b=T_{a,b}$.

\mk

Assume now that $m$ is odd, then according to Brady and McCammond (see~\cite{brady_mccammond_artinthree}), there is an interesting presentation of $I_2(m)$ given by $I_2(m) = \<a,b \st w_m(a,b) = w_m(b,a)\> = \<t,u \st t^m=u^2\>$, where $t=ab$ and $u=w_m(a,b)$, so the central quotient $G$ of $I_2(m)$ is isomorphic to $ \<t,u \st t^m=u^2\> / \<t^m=u^2\> \simeq \Z/m\Z \star \Z/2\Z$. Consider the action of $G$ on the Bass-Serre $(m,2)$-biregular tree $T$, and remove all degree $2$ vertices to obtain an action on the $m$-regular tree $T_{a,b}$. Note that since $b=uau^{-1}$, the axes $T_a,T_b$ of $a,b$ in $T_{a,b}$ intersect along the edge correponding to $\<u\>$. The two endpoints of this edge are $t_0$, corresponding to the subgroup $\<utu^{-1}\>=\<ba\>$, and $t_1$, corresponding to the subgroup $\<t\>=\<ab\>$. We have $a \cdot t_0 = t_1$ and $b \cdot t_1=t_0$. See Figure~\ref{fig:T3}.

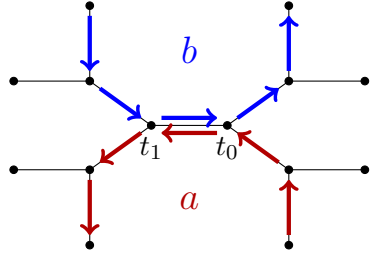
\begin{figure}
\begin{center}
\begin{tikzpicture}
\def \p {0.05}
\def \op {1}
\def \gris {black!10}
\draw[fill] (-0.5,0) circle (\p) node(1) {};
\draw[fill] (1)+(144:1) circle (\p) node(2) {};
\draw[fill] (1)+(-144:1) circle (\p) node(5) {};
\draw[fill] (2)+(90:1) circle (\p) node(4) {};
\draw[fill] (2)+(180:1) circle (\p) node(3) {};
\draw[fill] (5)+(-90:1) circle (\p) node(7) {};
\draw[fill] (5)+(-180:1) circle (\p) node(6) {};
\draw[fill] (0.5,0) circle (\p) node(0) {};
\draw[fill] (0)+(36:1) circle (\p) node(8) {};
\draw[fill] (0)+(-36:1) circle (\p) node(11) {};
\draw[fill] (8)+(0:1) circle (\p) node(9) {};
\draw[fill] (8)+(90:1) circle (\p) node(10) {};
\draw[fill] (11)+(0:1) circle (\p) node(12) {};
\draw[fill] (11)+(-90:1) circle (\p) node(13) {};
\node (0S) at (0.5,-0.1) {};
\node (0N) at (0.5,0.1) {};
\node (1S) at (-0.5,-0.1) {};
\node (1N) at (-0.5,0.1) {};

\draw[black] (7.center) -- (5.center);
\draw[black] (6.center) -- (5.center) -- (1.center) -- (2.center) -- (3.center);
\draw[black] (2.center) -- (4.center);
\draw[black] (1.center) -- (0.center) -- (8.center) -- (10.center) (8.center) -- (9.center) (0.center) -- (11.center) -- (12.center) (11.center) -- (13.center);

\draw [->,ultra thick,black!30!red] (13) edge (11) (11) edge (0) (0S) edge (1S) (1) edge (5) (5) edge (7);
\draw [->,ultra thick,blue] (4) edge (2) (2) edge (1) (1N) edge (0N) (0) edge (8) (8) edge (10);

\node (a) at (0,-1) {\color{black!30!red} \Large\bfseries $a$};
\node (a) at (0,1) {\color{blue} \Large\bfseries $b$};
\node (t0) at ([yshift=-0.3cm]0) {\bfseries $t_0$};
\node (t1) at ([yshift=-0.3cm]1) {\bfseries $t_1$};

\end{tikzpicture}
\end{center}
\caption{A part of the tree $T_{a,b}$ for $m_{a,b}=3$, with the axes of $a$ and $b$.}
\label{fig:T3}
\end{figure}

\mk

Assume finally that $m=2p$ is even, then according to Brady and McCammond (see~\cite{brady_mccammond_artinthree}), there is an interesting presentation of $I_2(m)$ given by $I_2(m) = \<a,b \st w_m(a,b) = w_m(b,a)\> = \<a,t \st at^p=t^pa\>$, where $t=ab$. In particular, $I_2(m)$ can be seen as the HNN extension of the group $\<t\> \simeq \Z$ with the subgroup $\<t^p\>$ and the identity map, with stable letter $a$.

Consider the action of $I_2(2p)$ on the Bass-Serre oriented $2p$-regular tree $T_{a,b}$. Let $t_0$ denote the vertex corresponding to the subgroup $\<a^{-1}ta\>$, it is fixed by $\<a^{-1}ta\>=\<ba\>$, and let $t_1 = a \cdot t_0$ denote the vertex corresponding to the subgroup $\<t\>$, it is fixed by $\<t\>=\<ab\>$. The axes $T_a,T_b$ of $a,b$ intersect in the edge $[t_0,t_1]$. See Figure~\ref{fig:T4}.

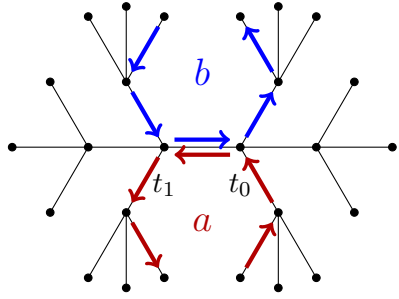
\begin{figure}
\begin{center}
\begin{tikzpicture}
\def \p {0.05}
\def \op {1}
\def \gris {black!10}
\draw[fill] (-0.5,0) circle (\p) node(1) {};
\draw[fill] (1)+(120:1) circle (\p) node(2) {};
\draw[fill] (1)+(-120:1) circle (\p) node(5) {};
\draw[fill] (1)+(180:1) circle (\p) node(23) {};
\draw[fill] (2)+(60:1) circle (\p) node(4) {};
\draw[fill] (2)+(90:1) circle (\p) node(3) {};
\draw[fill] (2)+(120:1) circle (\p) node(21) {};
\draw[fill] (5)+(-60:1) circle (\p) node(7) {};
\draw[fill] (5)+(-90:1) circle (\p) node(6) {};
\draw[fill] (5)+(-120:1) circle (\p) node(22) {};
\draw[fill] (23)+(120:1) circle (\p) node(26) {};
\draw[fill] (23)+(180:1) circle (\p) node(25) {};
\draw[fill] (23)+(-120:1) circle (\p) node(24) {};
\draw[fill] (0.5,0) circle (\p) node(0) {};
\draw[fill] (0)+(60:1) circle (\p) node(8) {};
\draw[fill] (0)+(-60:1) circle (\p) node(11) {};
\draw[fill] (0)+(0:1) circle (\p) node(33) {};
\draw[fill] (8)+(90:1) circle (\p) node(9) {};
\draw[fill] (8)+(120:1) circle (\p) node(10) {};
\draw[fill] (8)+(60:1) circle (\p) node(31) {};
\draw[fill] (11)+(-90:1) circle (\p) node(12) {};
\draw[fill] (11)+(-120:1) circle (\p) node(13) {};
\draw[fill] (11)+(-60:1) circle (\p) node(32) {};
\draw[fill] (33)+(60:1) circle (\p) node(36) {};
\draw[fill] (33)+(0:1) circle (\p) node(35) {};
\draw[fill] (33)+(-60:1) circle (\p) node(34) {};

\node (0S) at (0.5,-0.1) {};
\node (0N) at (0.5,0.1) {};
\node (1S) at (-0.5,-0.1) {};
\node (1N) at (-0.5,0.1) {};

\draw[black] (7.center) -- (5.center);
\draw[black] (6.center) -- (5.center) -- (1.center) -- (2.center) -- (3.center);
\draw[black] (2.center) -- (4.center);
\draw[black] (1.center) -- (0.center) -- (8.center) -- (10.center) (8.center) -- (9.center) (0.center) -- (11.center) -- (12.center) (11.center) -- (13.center);
\draw[black] (2.center) -- (21.center) (5.center) -- (22.center) (1.center) -- (23.center) -- (25.center) (24.center) -- (23.center) -- (26.center);
\draw[black] (8.center) -- (31.center) (11.center) -- (32.center) (0.center) -- (33.center) -- (35.center) (34.center) -- (33.center) -- (36.center);

\draw [->,ultra thick,black!30!red] (13) edge (11) (11) edge (0) (0S) edge (1S) (1) edge (5) (5) edge (7);
\draw [->,ultra thick,blue] (4) edge (2) (2) edge (1) (1N) edge (0N) (0) edge (8) (8) edge (10);

\node (a) at (0,-1) {\color{black!30!red} \Large\bfseries $a$};
\node (a) at (0,1) {\color{blue} \Large\bfseries $b$};
\node (t0) at ([yshift=-0.5cm]0) {\bfseries $t_0$};
\node (t1) at ([yshift=-0.5cm]1) {\bfseries $t_1$};

\end{tikzpicture}
\end{center}
\caption{A part of the tree $T_{a,b}$ for $m_{a,b}=4$, with the axes of $a$ and $b$.}
\label{fig:T4}
\end{figure}
\ep

\section{The gluing construction - cases (A) and (B)} \label{sec:AB_glue}

Consider a labeled graph $\Gamma$ such that:
\ben
\item[(A)] either $\Gamma$ has no cycle
\item[(B)] or $\Gamma$ is bipartite, and all its labels are at least $3$.
\een	

In both cases, $\Gamma$ is bipartite, and we will fix a coloring $\chi$ of vertices of $\Gamma$ in $\{0,1\}$ such that adjacent vertices have distinct colors.

If $\Gamma$ is not connected, then the Artin group $A(\Gamma)$ is the free product of the parabolic subgroups corresponding to the connected components of $\Gamma$. In order to prove Theorem~\ref{thm:main}, it is enough to consider the case where $\Gamma$ is connected.

For each $a \in S$, we will denote by $T_a \simeq \R$ the real line tiled by unit segments, with particular vertices $t_0=0$ and $t_1=1$. We also define the action of $\<a\>$ on $T_a$ where $a$ acts as a translation by $(-1)^{\chi(a)}$.

For each $a \in S$, let $X_a = T_a \times \Sigma_a \simeq \R^2$, with the product action of $a$, the product cubical structure, and with base vertex $p_0=(t_0,x_0)$ and its particular neighbour $p_1=(t_1,x_0)$. Let us denote the quotient $M_a = \<a\> \bs X_a$, with base vertex $q_0$, the image of $p_0$, and its particular neighbour $q_1$, image of $p_1$.

\blem
The action of $\<a\>$ on $X_a$ is free, so $M_a$ is a locally CAT(0) square complex with two hyperplanes, and the fundamental group $\pi_1(M_a,q_0)$ is naturally isomorphic to $\<a\>$.
\elem

\bp
Since the action of $\<a\>$ on the factor $T_a \simeq \R$ is free, we deduce that the action of $\<a\>$ on $X_a$ is free. Notice that $\<a\>$ acts transitively on the hyperplanes of $T_a$ and on the hyperplanes of $\Sigma_a$, so that $M_a$ has two hyperplanes.
\ep

Let $E$ denote the set of edges of $\Gamma$. For each edge $\{a,b\} \in E$, let $X_{a,b} = T_{a,b} \times \Sigma_{a,b}$, where $T_{a,b}$ denotes the tree described in Lemma~\ref{lem:tree_action_A} for the dihedral Artin group $\<a,b\>$, and $\Sigma_{a,b}$ denotes the Euclidean space described in Proposition~\ref{pro:euclidean_action} for the dihedral Artin group $\<a,b\>$. Note that $X_{a,b}$ is endowed with the product cubical structure, and the product action of $\<a,b\>$. It has a base vertex $p_0=(t_0,x_0)$ and its particular neighbour $p_1=(t_1,x_0)$. Let us denote the quotient $M_{a,b} = \<a,b\> \bs X_{a,b}$, with base vertex $q_0$, the image of $p_0$, and its particular neighbour $q_1$, image of $p_1$.

\blem
The action of $\<a,b\>$ on $X_{a,b}$ is free, so $M_{a,b}$ is a locally CAT(0) cube complex of dimension $m_{a,b}+1$, with two or three hyperplanes, and the fundamental group $\pi_1(M_{a,b},q_0)$ is naturally isomorphic to $\<a,b\>$.
\elem

\bp
Assume that $(t,x) \in X_{a,b}$ and $g \in \<a,b\>$ are such that $g \cdot (t,x)=(t,x)$. We will prove that $g=1$.

Since $\<a,b\>$ acts transitively on the vertices and on the edges of $T_{a,b}$, we may assume that either $t=t_0$ or $t \in (t_0,t_1)$. Without loss of generality, assume that $\chi(a)=0$ and $\chi(b)=1$.

If $t=t_0$, then $g \in \<ba\>$. According to Lemma~\ref{lem:action_free}, $g=1$.

If $t \in (t_0,t_1)$, then the stabilizer of $t$ is contained in $\<w_m(a,b)\>$, so $g \in \<w_m(a,b)\>$. According to Lemma~\ref{lem:action_separates}, $g=1$.

\mk

Notice that $\<a,b\>$ acts transivitely on hyperplanes of $T_{a,b}$. If $m_{a,b}$ is odd, then $\<a,b\>$ acts transitively on hyperplanes of $\Sigma_{a,b}$, and if $m_{a,b}$ is even, then $\<a,b\>$ has two orbits of hyperplanes in $\Sigma_{a,b}$. Therefore $M_{a,b}$ has two or three hyperplanes.
\ep

Let us try to give a simple picture of the space $M_{ab}$ when $m=3$. Let $P_{a,b} < \<a,b\>$ denotes the pure braid subgroup, i.e. the kernel of the canonical morphism $\<a,b \st aba=bab\> \ra W_{a,b} = \<a,b \st a^2=1,b^2=1,aba=bab\> \simeq \frak{S}_3$. Rather than $M_{a,b}$, it is probably simpler to describe the finite cover $\widehat{M_{ab}} = X_{a,b} / P_{a,b}$. Note that the center $Z_{a,b}$ of $P_{ab}$ is the infinite cyclic subgroup spanned by $(ab)^3$, and that $P_{a,b}/Z_{a,b}$ acts freely, properly and cocompactly on the $3$-regular tree $T_{a,b}$, with quotient the graph $\ov{T_{a,b}}$ with two vertices and three edges, see Figure~\ref{fig:Mab3}. So we can describe $\widehat{M_{ab}}$ as a fiber bundle over $\ov{T_{a,b}}$, with fiber $\Sigma_{a,b}/Z_{a,b} \simeq \R^2 \times \SS^1$, and it is homeormorphic to the trivial bundle $\R^2 \times \SS^1 \times \ov{T_{a,b}}$.

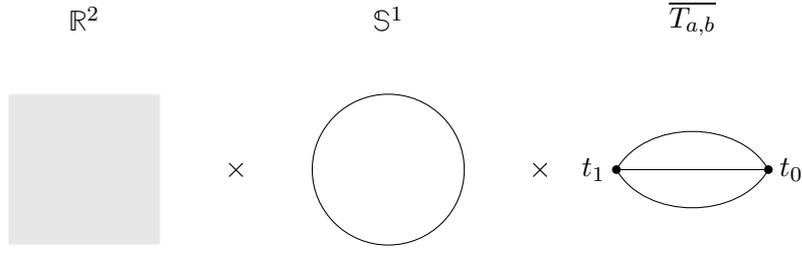
\begin{figure}
\begin{center}
\begin{tikzpicture}
\def \p {0.05}
\def \op {1}
\def \gris {black!10}
\draw[fill] (1,0) circle (\p) node(0) {};
\draw[fill] (-1,0) circle (\p) node(1) {};
\node (a) at (-9,-1) {};
\node (b) at (-7,-1) {};
\node (c) at (-7,1) {};
\node (d) at (-9,1) {};

\draw[black] (0.center) -- (1.center);
\path (0.center) edge [bend left=60] node {} (1.center);
\path (1.center) edge [bend left=60] node {} (0.center);
\draw[white,fill opacity=0.1,fill=black] (a.center) -- (b.center) -- (c.center) -- (d.center) -- (a.center);

\draw (-4,0) circle (1);

\node (t0) at ([xshift=0.3cm]0) {\bfseries $t_0$};
\node (t1) at ([xshift=-0.3cm]1) {\bfseries $t_1$};
\node (T) at (0,2) {$\ov{T_{a,b}}$};
\node (R) at (-8,2) {$\R^2$};
\node (S) at (-4,2) {$\SS^1$};
\node (X1) at (-2,0) {$\times$};
\node (X2) at (-6,0) {$\times$};

\end{tikzpicture}
\end{center}
\caption{A picture of the finite cover $\widehat{M_{ab}}$ of $M_{a,b}$ when $m=3$}
\label{fig:Mab3}
\end{figure}

\blem \label{lem:injection_a_ab_A}
For each $a \in S$ and $\{a,b\} \in E$, the $\<a\>$-equivariant embedding $X_a \ra X_{a,b}$ defines an embedding $\phi_{a,\{a,b\}} : M_a \ra M_{a,b}$ sending $q_0,q_1 \in M_a$ to $q_0,q_1 \in M_{a,b}$. Furthermore, $\phi_{a,\{a,b\}}(M_a)$ is locally convex in $M_{a,b}$.
\elem

\bp
Assume that $(t,x)$ is a vertex in $X_a$ and $g \in \<a,b\>$ are such that $g \cdot (t,x) \in X_a$, we will show that $g \cdot (t,x) \in \<a\> \cdot (t,x)$.

The action of $\<a\>$ is transitive on vertices of $T_a$, so we may assume that $g \cdot t=t=t_0$. Without loss of generality, assume that $\chi(a)=0$ and $\chi(b)=1$. Then $g \in \<ba\>$, and there exists $n \in \Z$ such that $g=(ba)^n$.

Since $(ba)^n \cdot \Sigma_a \cap \Sigma_a \neq \emptyset$, according to Lemma~\ref{lem:action_separates}, we deduce that $n=0$, so $g=1$.

Furthermore, note that the image of the natural embedding $X_a \ra X_{a,b}$ has convex image in the CAT(0) cube complex $X_{a,b}$. So its image $\phi_{a,\{a,b\}}(M_a)$ in the quotient $M_{a,b}$ by the free action of $\<a,b\>$ is locally convex.\ep

Now define
$$M = \left( \bigcup_{a \in S} M_a \cup \bigcup_{\{a,b\} \in E} M_{a,b}\right) / \sim,$$
where the identifications are given, for every $a \in S$ and $\{a,b\} \in E$, by $\phi_{a,\{a,b\}} : M_a \ra M_{a,b}$.

\mk

It is a cube complex, with a basepoint $q_0$, and a particular neighbour $q_1$. We will first prove that each $M_a$ and $M_{a,b}$ embeds in $M$, and then prove that the fundamental group of $M$ is the Artin group $A(\Gamma)$.

\blem \label{lem:intersection_3_A}
For each edge $\{a,b\} \in E$ with label $m_{a,b} \geq 3$, we have
$$\phi_{a,\{a,b\}}(M_a) \cap \phi_{b,\{a,b\}}(M_b) = [q_0,q_1].$$
More precisely, $\phi_{a,\{a,b\}}(q_0)=\phi_{b,\{a,b\}}(q_0)$ and $\phi_{a,\{a,b\}}(q_1)=\phi_{b,\{a,b\}}(q_1)$.
\elem

\bp
Assume that $(t,x) \in X_a$ and $g \in \<a,b\>$ are such that $g \cdot (t,x) \in X_b$, we will prove that $(t,x) \in \{p_0,p_1\}$. Without loss of generality, assume that $\chi(a)=0$ and $\chi(b)=1$. 

Since $\<a\>$ acts transitively on the vertices of $T_a$, we may assume that $t=t_0$. Since $\<b\>$ acts transitively on the vertices of $T_b$, we may also assume that $g \cdot t=g \cdot t_0=t_0$. Therefore $g \in \<ba\>$, and there exists $n \in \Z$ such that $g=(ba)^n$. Since $(ba)^n \cdot \Sigma_a \cap \Sigma_b \neq \emptyset$, according to Lemma~\ref{lem:action_separates}, we deduce that $n=0$ and $g=1$.
\ep

\blem \label{lem:intersection_2_A}
For each edge $\{a,b\} \in E$ with label $m_{a,b} =2$, we have
$$\phi_{a,\{a,b\}}(M_a) \cap \phi_{b,\{a,b\}}(M_b) = \phi_{a,\{a,b\}}(T_a \times \{x_0\}) = \phi_{b,\{a,b\}}(T_b \times \{x_0\}).$$
\elem

\bp
Notice that every vertex in $\phi_{a,\{a,b\}}(M_a)$ is in the image of $\phi_{a,\{a,b\}}(T_a \times \{x_0\})$, but not every edge.

Consider an edge $[(t,x),(t',x')] \in X_a$ and $g \in \<a,b\>$ such that $g \cdot [(t,x),(t',x')] \in X_b$. We will prove that $g=1$ and that $x=x'=x_0$. Note that edges in $\Sigma_a$ and edges in $\Sigma_b$ are not in the same $\<a,b\>$-orbit, so $x=x'$. 

\mk

Since $\<a\>$ acts transitively on the vertices of $T_a$, we may assume that $t=t_0$. Since $\<b\>$ acts transitively on the vertices of $T_b$, we may also assume that $g \cdot t_0=t_0$. Therefore $g \in \<ba\>$, and there exists $n \in \Z$ such that $g=(ba)^n$. Now $(ba)^n \cdot \Sigma_a \cap \Sigma_b \neq \emptyset$, so according to Lemma~\ref{lem:action_separates} we deduce that $n=0$ and $g=1$. Therefore $x \in \Sigma_a \cap \Sigma_b = \{x_0\}$.
\ep

\blem \label{lem:M_ab_injective}
For each $a \in S$, the natural map $M_a \ra M$ is injective. For each $\{a,b\} \in E$, the natural map $M_{a,b} \ra M$ is injective.
\elem

\bp
Note that, according to Lemma~\ref{lem:injection_a_ab_A}, for every $a \in S$, the edge $[q_0,q_1] \subset M_a$ injects in $M$. Therefore, in order to simplify notation, we will consider $[q_0,q_1]$ as an edge of $M$.

\mk

We will first prove the following result. Fix $a,b \in S$, with either $a=b$ or $\{a,b\} \in E$. Assume that there exists a sequence $a_0=a,a_1,\dots,a_n=b$ in $S$ and points $r_0 \in M_{a_0},r_1 \in M_{a_1},\dots,r_n \in M_{a_n}$ such that for every $0 \leq i \leq n-1$, we have $\{a_i,a_{i+1}\} \in E$ and $\phi_{a_i,\{a_i,a_{i+1}\}}(r_i) = \phi_{a_{i+1},\{a_i,a_{i+1}\}}(r_{i+1})$. We will prove that
\ben
\item $a=b$ and $r_0=r_n \in M_a=M_b$, or
\item $\{a,b\} \in E$, $m_{a,b} \geq 3$  and $r_0=r_n \in [q_0,q_1]$, or
\item $\{a,b\} \in E$, $m_{a,b}=2$ and $r_0=r_n \in T_a \times \{x_0\}=T_b \times \{x_0\}=T_{a,b} \times \{x_0\}$.
\een
To prove that, we consider each of the two possibilities for $\Gamma$.
\ben
\item[(A)] Consider the case where $\Gamma$ is a tree. Then there exists $0 \leq i \leq n-2$ such that $a_i=a_{i+2}$. Since $M_{a_i}$ embeds in $M_{a_i,a_{i+1}}$ according to Lemma~\ref{lem:injection_a_ab_A}, we deduce that $r_i=r_{i+2}$. So we can shorten the path from $r_0$ to $r_n$. By induction on $n$, we are reduced to the case $n \leq 1$. If $n=0$ then $a=b$ and $r_0=r_n$. If $n=1$ and $m_{a,b} \geq 3$, then according to Lemma~\ref{lem:intersection_3_A} we deduce that $r_0=r_n \in [q_0,q_1]$. If $n=1$ and $m_{a,b}=2$, then according to Lemma~\ref{lem:intersection_2_A} we deduce that  $r_0=r_n \in T_a \times \{x_0\}=T_b \times \{x_0\}=T_{a,b} \times \{x_0\}$.
\item[(B)] Consider the case where $\Gamma$ is bipartite, and all labels are at least $3$. If $n=0$, then $a=b$ and $r_0=r_n$. If $n \geq 1$, then for each $0 \leq i \leq n-1$, as we have $\phi_{a_i,\{a_i,a_{i+1}\}}(r_i) = \phi_{a_{i+1},\{a_i,a_{i+1}\}}(r_{i+1})$, we deduce according to Lemma~\ref{lem:intersection_3_A} that $r_i,r_{i+1} \in [q_0,q_1]$. Therefore $r_0,r_n \in [q_0,q_1]$. According to Lemma~\ref{lem:injection_a_ab_A}, the edge $[q_0,q_1]$ injects in $M$, so $r_0=r_n$.
\een

\mk

For each $a \in S$, applying the previous result to $a=b \in S$ proves that $M_a$ embeds in $M$.

For each $\{a,b\} \in E$, applying the previous result to $a,b$ proves that $M_a \cup M_b$ injects in $M$. Since the subset $M_{a,b} \bs (\phi_{a,\{a,b\}}(M_a) \cup \phi_{b,\{a,b\}}(M_b))$ injects in $M$, we conclude that $M_{a,b}$ injects in $M$.
\ep

In order to simplify notation, we will therefore identify each $M_a$, for $a \in S$, and each $M_{a,b}$, for $\{a,b\} \in E$, with their images in $M$.

\blem
The fundamental group $\pi_1(M,q_0)$ is naturally isomorphic to $A=A(\Gamma)$.
\elem

\bp
We apply the Van Kampen Theorem to the subsets $M_{a,b}$, for $\{a,b\} \in E$, each containing the basepoint $q_0$. According to the proof of Lemma~\ref{lem:M_ab_injective}, for any two distinct edges $\{a,b\},\{c,d\}$ in $E$, the intersection $M_{a,b} \cap M_{c,d}$ can be $M_s$ (for some $s \in S$), or the image of $T_s \times \{x_0\}$ in $M$ (for some $s \in S$), or $[q_0,q_1]$. As a consequence, for any three distinct edges $\{a,b\},\{c,d\},\{e,f\}$ in $E$, the triple intersection $M_{a,b} \cap M_{c,d} \cap M_{e,f}$ is of the same form. We deduce that pairwise and triple intersections of subsets $M_{a,b}$, for $\{a,b\} \in E$, are connected. We conclude that the fundamental group $\pi_1(M,q_0)$ is naturally isomorphic to the Artin group $A=A(\Gamma)$.
\ep

\blem
In the case (A) where $\Gamma$ is a tree, the cube complex $M$ is locally CAT(0).
\elem

\bp
We will prove it by induction on the number of vertices of $\Gamma$. If $\Gamma$ is a single edge $\{a,b\}$, then $M=M_{a,b}$ is locally CAT(0). Otherwise, consider an edge $\{a,b\} \in E$ containing a leaf $a \in S$ of $\Gamma$. Let $\Gamma'$ denote the subtree obtained by removing $a$ from $\Gamma$. Let $M'$ denote the complex associated to $\Gamma'$. We have $M = M' \cup M_{a,b} \bs \sim$, where the identification is given by the two embeddings of $M_b$ in $M'$ and $M_{a,b}$. We know that $M_{a,b}$ is locally CAT(0) and $M'$ is locally CAT(0) by induction. Furthermore, $M_b$ is locally convex in both $M_{a,b}$ and $M'$. Therefore, $M$ is locally CAT(0).
\ep

\blem
In the case (B) where $\Gamma$ is bipartite and has no label $2$, the cube complex $M$ is locally CAT(0).
\elem

\bp
Fix a vertex $q \in M$, we will prove that $M$ is locally CAT(0) at $q$.

\mk

Assume first that there exists $\{a,b\} \in E$ such that $q \in M_{a,b} \bs \left(M_a \cup M_b\right)$. Then $M_{a,b}$ is a neighbourhood of $q$ in $M$. Since $M_{a,b}$ is locally CAT(0), $M$ is locally CAT(0) at $q$.

Assume now that there exists $a \in S$ such that $q \in M_a \bs [q_0,q_1]$. Then $M'_a = \bigcup_{\{a,b\} \in E} M_{a,b}$ is a neighbourhood of $q$ in $M$. Since $M'_a$ is the gluing of locally CAT(0) cube complexes $M_{a,b}$, for $\{a,b\} \in E$, along the common locally convex subspace $M_a$, we deduce that $M'_a$ is a locally CAT(0) cube complex. Therefore, $M$ is locally CAT(0) at $q$.

Assume now that $q=q_0$. We will prove that the link of $q_0$ is a flag simplicial complex. Assume that $Q_1,Q_2,Q_3$ are three cubes of $M$ containing $q_0$ such that each pairwise intersection $Q_i \cap Q_j$ has codimension $1$ in $Q_i$ and $Q_j$, and $Q_1 \cap Q_2 \cap Q_3$ has codimension $2$ in $Q_1$, $Q_2$ and $Q_3$. We will prove that there exists a cube $Q$ of $M$ containing $Q_1$, $Q_2$ and $Q_3$ with codimension $1$. For each $i \in \{1,2,3\}$, let $\{a_i,b_i\} \in E$ such that $Q_i \subset M_{a_i,b_i}$. 

\mk

{\bf Claim: } If there exist three edges $\{a,b\},\{a,c\},\{a,d\} \in E$ such that $Q_1,Q_2,Q_3 \subset M_{a,b} \cup M_{a,c} \cup M_{a,d}$, then such a cube $Q$ exists. Indeed note that $M_{a,b} \cup M_{a,c} \cup M_{a,d}$ is the union of three locally CAT(0) cube complexes along the common locally convex subspace $M_a$, hence it is locally CAT(0). Therefore there exists a cube $Q$ in $M_{a,b} \cup M_{a,c} \cup M_{a,d}$ containing each of $Q_1,Q_2,Q_3$ with codimension $1$ in this situation.

\mk

\bit
\item Assume first that for each $i \neq j$, we have $Q_i \cap Q_j \neq [q_0,q_1]$. Then for each $i \neq j$, $M_{a_i,b_i}$ and $M_{a_j,b_j}$ intersect outside of $[q_0,q_1]$, so the edges $\{a_i,b_i\}$ and $\{a_j,b_j\}$ intersect in $\Gamma$. Since $\Gamma$ has no triangles, there exist $\{a,b\},\{a,c\},\{a,d\} \in E$ such that $Q_1,Q_2,Q_3 \subset M_{a,b} \cup M_{a,c} \cup M_{a,d}$. According to the claim, a cube $Q$ as required exists.
\item Assume now that, for instance, we have $Q_1 \cap Q_3 = [q_0,q_1]$. Then, as $Q_1 \cap Q_2 \neq [q_0,q_1]$ and $Q_2 \cap Q_3 \neq [q_0,q_1]$, by the previous argument we know that the edges $\{a_1,b_1\}$ and $\{a_2,b_2\}$ in $\Gamma$, and also the edges $\{a_2,b_2\}$ and $\{a_3,b_3\}$ intersect in $\Gamma$. For instance, we can assume for instance that $b_1=a_2$. Since $Q_1 \cap Q_3 = [q_0,q_1]$, we know that $Q_1,Q_2,Q_3$ are squares. So $Q_1$ is the square at $q_0$ spanned by the edges $Q_1 \cap Q_3=[q_0,q_1]$ and $Q_1 \cap Q_2 \subset M_{a_1,b_1} \cap M_{a_2,b_2} = M_{a_2}$. So we deduce that $Q_1 \subset M_{a_2}$.  Hence we deduce that $Q_1,Q_2,Q_3$ are contained in $M_{a_2,b_2} \cup M_{a_3,b_3}$. Since the edges $\{a_2,b_2\}$ and $\{a_3,b_3\}$ intersect in $\Gamma$, by the claim we know that a cube $Q$ as desired exists.
\eit

Assume finally that $q=q_1$. This situation is entirely similar to the previous one $q=q_0$.
\ep

\section{The tree actions, case (C)} \label{sec:C_tree}

We will now describe an action by isometries of each dihedral Artin group on a tree-like CAT(0) square complex, that will be used in the case (C).

\blem \label{lem:tree_action_B}
For every $m \neq 3$, there exists a cubical action of the dihedral Artin group $I_2(m)=\<a,b \st w_m(a,b)=w_m(b,a)\>$ on a CAT(0) square complex $T_{a,b}$ such that the following hold.
\ben
\item The elements $a,b$ act as translations on $T_{a,b}$, with combinatorial displacement $2$, and with axes combinatorial lines $T_a,T_b$ such that $T_a \cap T_b$ is a single vertex $t_0$.
\item If $m \geq 4$, the stabilizer of $t_0$ is $\<w_m(a,b)\>$. If $m=2$, the stabilizer of $t_0$ is trivial.
\een
\elem

\bp

Assume first that $m=2p+1 \geq 5$ is odd, then according to Brady and McCammond (see~\cite{brady_mccammond_artinthree}), there is an interesting presentation of $I_2(m)$ given by $I_2(m) = \<a,b \st w_m(a,b) = w_m(b,a)\> = \<t,u \st t^m=u^2\>$, where $t=ab$ and $u=w_m(a,b)$, so the central quotient $G$ of $I_2(m)$ is isomorphic to $ \<t,u \st t^m=u^2\> / \<t^m=u^2\> \simeq \Z/m\Z \star \Z/2\Z$. Consider the action of $G$ on the Bass-Serre $(m,2)$-biregular tree $T$. Consider the square complex $T_{a,b}$ obtained from $T$ by replacing the star of each vertex with valency $m$ by a regular $m$-gon tessalated by $m$ squares, where $t$ acts on the base $m$-gon $P_0$ by a rotation of angle $\frac{4\pi}{m}$. Note that $a=t^{-p}u$ and $b=ut^{-p}$, and $t^p$ acts on the base $m$-gon by a rotation of angle $\frac{4p\pi}{m}=\frac{-2\pi}{m}$. This way, the axes of $a$ and $b$ acting on $T_{a,b}$ intersect the boundary of the $m$-gon $P$ in consecutive sides. Let $t_0 \in T_{a,b}$ denote the intersection of the axes of $a$ and $b$, it is also the unique vertex fixed by $u=w_m(a,b)$ (see Figure~\ref{fig:T5}).

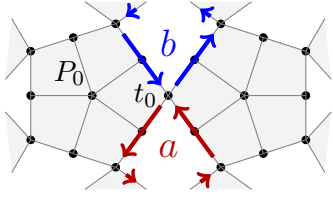
\begin{figure}
\begin{center}
\begin{tikzpicture}
\def \p {0.05}
\def \op {0.05}
\def \l {0.809}
\def \gris {black}
\def \noir {black!50}
\draw[fill] (-1,0) circle (\p) node(P1) {};
\draw[fill] (-1,0)+(0:1) circle (\p) node(1) {};
\draw[fill] (-1,0)+(36:\l) circle (\p) node(1m2) {};
\draw[fill] (-1,0)+(72:1) circle (\p) node(2){} ;
\draw[fill] (-1,0)+(108:\l) circle (\p) node(2m3) {};
\draw[fill] (-1,0)+(144:1) circle (\p) node(3) {};
\draw[fill] (-1,0)+(180:\l) circle (\p) node(3m4) {};
\draw[fill] (-1,0)+(-144:1) circle (\p) node(4) {};
\draw[fill] (-1,0)+(-108:\l) circle (\p) node(4m5) {};
\draw[fill] (-1,0)+(-72:1) circle (\p) node(5) {};
\draw[fill] (-1,0)+(-36:\l) circle (\p) node(5m1) {};
\draw[fill] (1,0) circle (\p) node(P2) {};
\draw[fill] (1,0)+(180+36:\l) circle (\p) node(1m6) {};
\draw[fill] (1,0)+(180+72:1) circle (\p) node(6) {};
\draw[fill] (1,0)+(180+108:\l) circle (\p) node(6m7) {};
\draw[fill] (1,0)+(180+144:1) circle (\p) node(7) {};
\draw[fill] (1,0)+(0:\l) circle (\p) node(7m8) {};
\draw[fill] (1,0)+(180-144:1) circle (\p) node(8) {};
\draw[fill] (1,0)+(180-108:\l) circle (\p) node(8m9) {};
\draw[fill] (1,0)+(180-72:1) circle (\p) node(9) {};
\draw[fill] (1,0)+(180-36:\l) circle (\p) node(9m1) {};

\node (31) at ([shift=(120:0.5)]3) {};
\node (32) at ([shift=(120+108:0.5)]3) {};
\node (41) at ([shift=(-120:0.5)]4) {};
\node (42) at ([shift=(-120-108:0.5)]4) {};
\draw[\noir,fill opacity=\op,fill=\gris] (31.center) -- (3.center) -- (32.center);
\draw[\noir,fill opacity=\op,fill=\gris] (41.center) -- (4.center) -- (42.center);

\node (21) at ([shift=(36:0.5)]2) {};
\node (22) at ([shift=(36+108:0.5)]2) {};
\node (51) at ([shift=(-36:0.5)]5) {};
\node (52) at ([shift=(-36-108:0.5)]5) {};
\draw[\noir,fill opacity=\op,fill=\gris] (21.center) -- (2.center) -- (22.center);
\draw[\noir,fill opacity=\op,fill=\gris] (51.center) -- (5.center) -- (52.center);

\node (81) at ([shift=(60:0.5)]8) {};
\node (82) at ([shift=(60-108:0.5)]8) {};
\node (71) at ([shift=(-60:0.5)]7) {};
\node (72) at ([shift=(-60+108:0.5)]7) {};
\draw[\noir,fill opacity=\op,fill=\gris] (81.center) -- (8.center) -- (82.center);
\draw[\noir,fill opacity=\op,fill=\gris] (71.center) -- (7.center) -- (72.center);

\node (61) at ([shift=(-36:0.5)]6) {};
\node (62) at ([shift=(-36-108:0.5)]6) {};
\node (91) at ([shift=(36:0.5)]9) {};
\node (92) at ([shift=(36+108:0.5)]9) {};
\draw[\noir,fill opacity=\op,fill=\gris] (61.center) -- (6.center) -- (62.center);
\draw[\noir,fill opacity=\op,fill=\gris] (91.center) -- (9.center) -- (92.center);

\draw[\noir,fill opacity=\op,fill=\gris] (1.center) -- (1m2.center) -- (2.center) -- (2m3.center) -- (3.center) -- (3m4.center) -- (4.center) -- (4m5.center) -- (5.center) -- (5m1.center) -- (1.center);
\draw[\noir] (P1.center) -- (1m2.center) (P1.center) -- (2m3.center) (P1.center) -- (3m4.center) (P1.center) -- (4m5.center) (P1.center) -- (5m1.center);
\draw[\noir,fill opacity=\op,fill=\gris] (1.center) -- (9m1.center) -- (9.center) -- (8m9.center) -- (8.center) -- (7m8.center) -- (7.center) -- (6m7.center) -- (6.center) -- (1m6.center) -- (1.center);
\draw[\noir] (P2.center) -- (1m6.center) (P2.center) -- (6m7.center) (P2.center) -- (7m8.center) (P2.center) -- (8m9.center) (P2.center) -- (9m1.center);

\draw [->,ultra thick,black!30!red] (62) edge (6) (6) edge (1) (1) edge (5) (5) edge (51);
\draw [<-,ultra thick,blue] (92) edge (9) (9) edge (1) (1) edge (2) (2) edge (21);

\node (P) at ([xshift=-1.3cm,yshift=0.3cm]1) {\bfseries $P_0$};
\node (e) at ([xshift=-0.3cm]1) {\bfseries $t_0$};
\node (a) at ([yshift=-0.7cm]1) {\color{black!30!red} \Large\bfseries $a$};
\node (b) at ([yshift=+0.7cm]1) {\color{blue} \Large\bfseries $b$};

\end{tikzpicture}
\end{center}
\caption{A part of the complex $T_{a,b}$ for $m_{a,b}=5$, with the axes of $a$ and $b$.}
\label{fig:T5}
\end{figure}

\mk

Assume now that $m=2p \geq 4$ is even, then according to Brady and McCammond (see~\cite{brady_mccammond_artinthree}), there is an interesting presentation of $I_2(m)$ given by $I_2(m) = \<a,b \st w_m(a,b) = w_m(b,a)\> = \<a,t \st at^p=t^pa\>$, where $t=ab$. In particular, $I_2(m)$ can be seen as the HNN extension of the group $\<t\> \simeq \Z$ with the subgroup $\<t^p\>$ and the identity map, with stable letter $a$.

Consider the action of $I_2(2p)$ on the Bass-Serre oriented $2p$-regular tree $T$. Let $T'$ denote the barycentric subdivision of $T$, it is an oriented $(2p,2)$-biregular tree. Consider the square complex $T_{a,b}$ obtained from $T'$ by replacing the star of each vertex with degree $2p$ by a regular $2p$-gon tessalated by $2p$ squares, such that $t$ acts on the base $2p$-gon $P_0$ by a rotation of angle $\frac{4\pi}{2p}$.

Since $b=a^{-1}t$, the axes of $a$ and $b$ acting on $T_{a,b}$ intersect the boundary of the $2p$-gon $P$ in consecutive sides. Let $t_0 \in T_{a,b}$ denote the intersection of the axes of $a$ and $b$ (see Figure~\ref{fig:T6}).

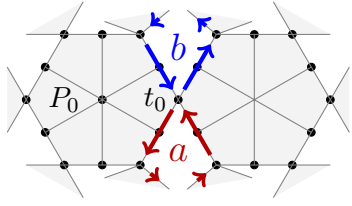
\begin{figure}
\begin{center}
\begin{tikzpicture}
\def \p {0.05}
\def \op {0.05}
\def \l {0.866}
\def \gris {black}
\def \noir {black!50}
\draw[fill] (-1,0) circle (\p) node(P1) {};
\draw[fill] (-1,0)+(0:1) circle (\p) node(1) {};
\draw[fill] (-1,0)+(30:\l) circle (\p) node(1m2) {};
\draw[fill] (-1,0)+(60:1) circle (\p) node(2){} ;
\draw[fill] (-1,0)+(90:\l) circle (\p) node(2m3) {};
\draw[fill] (-1,0)+(120:1) circle (\p) node(3) {};
\draw[fill] (-1,0)+(150:\l) circle (\p) node(3m4) {};
\draw[fill] (-1,0)+(180:1) circle (\p) node(4) {};
\draw[fill] (-1,0)+(-150:\l) circle (\p) node(4m5) {};
\draw[fill] (-1,0)+(-120:1) circle (\p) node(5) {};
\draw[fill] (-1,0)+(-90:\l) circle (\p) node(5m6) {};
\draw[fill] (-1,0)+(-60:1) circle (\p) node(6) {};
\draw[fill] (-1,0)+(-30:\l) circle (\p) node(6m1) {};
\draw[fill] (1,0)+(180+30:\l) circle (\p) node(1m7) {};
\draw[fill] (1,0)+(180+60:1) circle (\p) node(7) {};
\draw[fill] (1,0)+(180+90:\l) circle (\p) node(7m8) {};
\draw[fill] (1,0)+(180+120:1) circle (\p) node(8) {};
\draw[fill] (1,0)+(180+150:\l) circle (\p) node(8m9) {};
\draw[fill] (1,0)+(180+180:1) circle (\p) node(9) {};
\draw[fill] (1,0)+(180-150:\l) circle (\p) node(9m10) {};
\draw[fill] (1,0)+(180-120:1) circle (\p) node(10) {};
\draw[fill] (1,0)+(180-90:\l) circle (\p) node(10m11) {};
\draw[fill] (1,0)+(180-60:1) circle (\p) node(11) {};
\draw[fill] (1,0)+(180-30:\l) circle (\p) node(11m1) {};

\draw[\noir,fill opacity=\op,fill=\gris] (1.center) -- (1m2.center) -- (2.center) -- (2m3.center) -- (3.center) -- (3m4.center) -- (4.center) -- (4m5.center) -- (5.center) -- (5m6.center) -- (6.center) -- (6m1.center) -- (1.center);
\draw[\noir] (P1.center) -- (1m2.center) (P1.center) -- (2m3.center) (P1.center) -- (3m4.center) (P1.center) -- (4m5.center) (P1.center) -- (5m6.center) (P1.center) -- (6m1.center);
\draw[\noir,fill opacity=\op,fill=\gris] (1.center) -- (11m1.center) -- (11.center) -- (10m11.center) -- (10.center) -- (9m10.center) -- (9.center) -- (8m9.center) -- (8.center) -- (7m8.center) -- (7.center) -- (1m7.center) -- (1.center);
\draw[\noir] (P2.center) -- (1m7.center) (P2.center) -- (7m8.center) (P2.center) -- (8m9.center) (P2.center) -- (9m10.center) (P2.center) -- (10m11.center) (P2.center) -- (11m1.center);

\node (2a) at ([shift=(40:0.5)]2) {};
\node (2b) at ([shift=(160:0.5)]2) {};
\draw[\noir,fill opacity=\op,fill=\gris] (2a.center) -- (2.center) -- (2b.center);
\node (3a) at ([shift=(60:0.5)]3) {};
\node (3b) at ([shift=(180:0.5)]3) {};
\draw[\noir,fill opacity=\op,fill=\gris] (3a.center) -- (3.center) -- (3b.center);
\node (4a) at ([shift=(120:0.5)]4) {};
\node (4b) at ([shift=(-120:0.5)]4) {};
\draw[\noir,fill opacity=\op,fill=\gris] (4a.center) -- (4.center) -- (4b.center);
\node (6a) at ([shift=(-40:0.5)]6) {};
\node (6b) at ([shift=(-160:0.5)]6) {};
\draw[\noir,fill opacity=\op,fill=\gris] (6a.center) -- (6.center) -- (6b.center);
\node (5a) at ([shift=(-60:0.5)]5) {};
\node (5b) at ([shift=(-180:0.5)]5) {};
\draw[\noir,fill opacity=\op,fill=\gris] (5a.center) -- (5.center) -- (5b.center);

\node (11a) at ([shift=(180-40:0.5)]11) {};
\node (11b) at ([shift=(180-160:0.5)]11) {};
\draw[\noir,fill opacity=\op,fill=\gris] (11a.center) -- (11.center) -- (11b.center);
\node (10a) at ([shift=(180-60:0.5)]10) {};
\node (10b) at ([shift=(180-180:0.5)]10) {};
\draw[\noir,fill opacity=\op,fill=\gris] (10a.center) -- (10.center) -- (10b.center);
\node (9a) at ([shift=(180-120:0.5)]9) {};
\node (9b) at ([shift=(180+120:0.5)]9) {};
\draw[\noir,fill opacity=\op,fill=\gris] (9a.center) -- (9.center) -- (9b.center);
\node (8a) at ([shift=(-120:0.5)]8) {};
\node (8b) at ([shift=(0:0.5)]8) {};
\draw[\noir,fill opacity=\op,fill=\gris] (8a.center) -- (8.center) -- (8b.center);
\node (7a) at ([shift=(-140:0.5)]7) {};
\node (7b) at ([shift=(-20:0.5)]7) {};
\draw[\noir,fill opacity=\op,fill=\gris] (7a.center) -- (7.center) -- (7b.center);

\draw [->,ultra thick,black!30!red] (7a) edge (7) (7) edge (1) (1) edge (6) (6) edge (6a);
\draw [->,ultra thick,blue] (2a) edge (2) (2) edge (1) (1) edge (11) (11) edge (11a);

\node (P) at ([xshift=-1.5cm]1) {\bfseries $P_0$};
\node (e) at ([xshift=-0.3cm]1) {\bfseries $t_0$};
\node (a) at ([yshift=-0.7cm]1) {\color{black!30!red} \Large\bfseries $a$};
\node (b) at ([yshift=+0.7cm]1) {\color{blue} \Large\bfseries $b$};

\end{tikzpicture}
\end{center}
\caption{A part of the complex $T_{a,b}$ for $m_{a,b}=6$, with the axes of $a$ and $b$.}
\label{fig:T6}
\end{figure}

When $m=2$, i.e. for the abelian dihedral group $I_2(2)=\<a,b \st ab=ba\> \simeq \Z^2$, let $T_{a,b}=\R^2$ with the usual square tiling, where $a$ acts as a translation of $(2,0)$, and $b$ acts as a translation of $(0,2)$. Then $T_a=\R \times \{0\}$ and $T_b = \{0\} \times \R$ intersect in $t_0=(0,0)$.
\ep

\section{The gluing construction, case (C)} \label{sec:C_glue}

Fix a triangle-free graph $\Gamma$, with each edge being labeled by an integer equal to $2$ or at least $4$.

For each $a \in S$, let $X_a = T_a \times \Sigma_a \simeq \R^2$, with the product action of $a$, the product cubical structure, and with base vertex $p_0=(t_0,x_0)$. Let us denote the quotient $M_a = \<a\> \bs X_a$, with base vertex $q_0$, the image of $p_0$

\blem
The action of $\<a\>$ on $X_a$ is free, so $M_a$ is locally CAT(0) square complex with three hyperplanes, and the fundamental group $\pi_1(M_a,q_0)$ is naturally isomorphic to $\<a\>$.
\elem

\bp
Since the action of $\<a\>$ on the factor $T_a \simeq \R$ is free, we deduce that the action of $\<a\>$ on $X_a$ is free. Notice that $\<a\>$ has two orbits of hyperplanes in $T_a$ and acts transitively on the hyperplanes of $\Sigma_a$, so that $M_a$ has three hyperplanes.
\ep

Let $E$ denote the set of edges of $\Gamma$. For each edge $\{a,b\} \in E$, let $X_{a,b} = T_{a,b} \times \Sigma_{a,b}$, where $T_{a,b}$ denotes the square complex described in Lemma~\ref{lem:tree_action_B} for the dihedral Artin group $\<a,b\>$, and $\Sigma_{a,b}$ denotes the Euclidean space described in Proposition~\ref{pro:euclidean_action} for the dihedral Artin group $\<a,b\>$. Note that $X_{a,b}$ is endowed with the product cubical structure, and the product action of $\<a,b\>$. It has a base vertex $p_0=(t_0,x_0)$. Let us denote the quotient $M_{a,b} = \<a,b\> \bs X_{a,b}$, with base vertex $q_0$, the image of $p_0$.

\blem
The action of $\<a,b\>$ on $X_{a,b}$ is free, so $M_{a,b}$ is a locally CAT(0) cube complex of dimension $m_{a,b}+2$, with two or four hyperplanes, and the fundamental group $\pi_1(M_{a,b},q_0)$ is naturally isomorphic to $\<a,b\>$.
\elem

\bp
Assume that $(t,x) \in X_{a,b}$ and $g \in \<a,b\>$ are such that $g \cdot (t,x)=(t,x)$. We will prove that $g=1$.

\mk

Let $t$ denote the center of a polygon of $T_{a,b}$. Note thate $\<a,b\>$ acts transitively on polygons of $T_{a,b}$. Therefore we may assume that $t$ is the center of $b \cdot P_0$, and so $g \in b\<ab\>b^{-1}=\<ba\>$. According to Lemma~\ref{lem:action_free}, $g=1$.

\mk

If $t$ is in the orbit of $t_0$, we may assume that $t=t_0$, and so $g \in \<w_m(a,b)\>$.

If $t$ is not the center of a polygon nor in the orbit of $t_0$, then $g \in \<w_m(a,b)^2\>$ (if $m=m_{a,b}$ is odd) or $g \in \<w_m(a,b)\>$ (if $m=m_{a,b}$ is even).

In these last two cases, we have $g \in \<w_m(a,b)\>$. According to Lemma~\ref{lem:action_free}, $g=1$.

\mk

Notice that if $m$ is odd, then $\<a,b\>$ acts transivitely on hyperplanes of $T_{a,b}$, and if $m$ is even, then $\<a,b\>$ has two orbits of hyperplanes in $T_{a,b}$. If $m_{a,b}$ is odd, then $\<a,b\>$ acts transitively on hyperplanes of $\Sigma_{a,b}$, and if $m_{a,b}$ is even, then $\<a,b\>$ has two orbits of hyperplanes in $\Sigma_{a,b}$. Therefore $M_{a,b}$ has two or four hyperplanes.
\ep

\blem \label{lem:injection_a_ab_B}
For each $a \in S$ and $\{a,b\} \in E$, the $\<a\>$-equivariant embedding $X_a \ra X_{a,b}$ defines an embedding $\phi_{a,\{a,b\}} : M_a \ra M_{a,b}$ sending $q_0 \in M_a$ to $q_0 \in M_{a,b}$. Furthermore, $\phi_{a,\{a,b\}}(M_a)$ is locally convex in $M_{a,b}$.
\elem

\bp
Assume that $(t,x)$ is a vertex in $X_a$ and $g \in \<a,b\>$ are such that $g \cdot (t,x) \in X_a$, we will show that $g \cdot (t,x) \in \<a\> \cdot (t,x)$.

The action of $\<a\>$ has two orbits of vertices on $T_a$, so we may assume that $g \cdot t=t=t_0$ or that $g \cdot t=t$ is the common neighbour $t_a$ to $t_0$ and $a \cdot t_0$.
\ben
\item If $g \cdot t_0=t_0$, then $g \in \<w_m(a,b)\>$.
\item If $g \cdot t_a=t_a$ and $m$ is odd, then $g \in \<w_m(a,b)^2\>$.
\item If $g \cdot t_a=t_a$ and $m \geq 4$ is even, then $g \in \<w_m(a,b)\>$.
\item If $g \cdot t_a=t_a$ and $m=2$, then $g=1$.
\een
In every situation, we have $g \in \<w_m(a,b)\>$, and there exists $n \in \Z$ such that $g=(w_m(a,b))^n$. According to Lemma~\ref{lem:action_separates}, since $w_m(a,b)^n \cdot \Sigma_a \cap \Sigma_a \neq \emptyset$, we have $n=0$ and $g=1$.

Furthermore, note that the image of $X_a $ under the natural embedding in $X_{a,b}$ is convex. So its image $\phi_{a,\{a,b\}}(M_a)$ in the quotient $M_{a,b}$ by the free action of $\<a,b\>$ is locally convex.\ep

Now define
$$M = \left( \bigcup_{a \in S} M_a \cup \bigcup_{\{a,b\} \in E} M_{a,b}\right) / \sim,$$
where the identifications are given, for every $a \in S$ and $\{a,b\} \in E$, by $\phi_{a,\{a,b\}} : M_a \ra M_{a,b}$.

\mk

The space $M$ is a cube complex, with a basepoint $q_0$. We will first prove that each $M_a$ and $M_{a,b}$ embeds in $M$, and then prove that the fundamental group of $M$ is the Artin group $A(\Gamma)$.

\blem \label{intersection_B}
For each edge $\{a,b\} \in E$, we have
$$\phi_{a,\{a,b\}}(M_a) \cap \phi_{b,\{a,b\}}(M_b) = \{q_0\}.$$
\elem

\bp
Assume that $(t,x) \in X_a$ and $g \in \<a,b\>$ are such that $g \cdot (t,x) \in X_b$, we will prove that $(t,x) = (t_0,x_0)$.

Since $\<a\>$ has two orbits of vertices on $T_a$, we may assume that $t=t_0$ or $t=t_a$, the common neighbour of $t_0$ and $a \cdot t_0$. Since $\<b\>$ has two orbits of vertices on $T_b$, we may assume that $g \cdot t=t_0$ or $g \cdot t=t_b$, the common neighbour of $t_0$ and $b \cdot t_0$.

\bit
\item Assume first that $g \cdot t=t=t_0$. Then $g \in \<w_m(a,b)\>$. 
\item Assume now that $t=t_a$ and $g \cdot t_a = t_b$. Since $t_a,t_b$ are in the same $\<a,b\>$-orbit in $T_{a,b}$, this only occurs when $m$ is odd. Then $t_b = w_m(a,b) \cdot t_a$, and the stabilizer of $t_a$ is $\<w_m(a,b)^2\>$, so $g \in w_m(a,b) \<w_m(a,b)^2\> \subset \<w_m(a,b)\>$.
\eit

In each case, we deduce that $g \in \<w_m(a,b)\>$. Since $w_m(a,b)^n \cdot \Sigma_a \cap \Sigma_b \neq \emptyset$, according to Lemma~\ref{lem:action_separates} we have $n=0$ and $g=1$.
\ep

\blem
For each $a \in S$, the natural map $M_a \ra M$ is injective. For each $\{a,b\} \in E$, the natural map $M_{a,b} \ra M$ is injective.
\elem

\bp
This is a direct consequence of Lemma~\ref{intersection_B}.
\ep

In order to simplify notation, we will therefore identify each $M_a$, for $a \in S$, and each $M_{a,b}$, for $\{a,b\} \in E$, with their images in $M$.

\blem
The fundamental group $\pi_1(M,q_0)$ is naturally isomorphic to $A=A(\Gamma)$.
\elem

\bp
We apply the Van Kampen Theorem to the subsets $M_{a,b}$, for $\{a,b\} \in E$, each containing the basepoint $q_0$. Since pairwise and triple intersections are connected, we conclude that the fundamental group $\pi_1(M,q_0)$ is naturally isomorphic to the Artin group $A=A(\Gamma)$.
\ep

\blem
The cube complex $M$ is locally CAT(0).
\elem

\bp
Fix a vertex $q \in M$, we will prove that $M$ is locally CAT(0) at $q$.

\mk

Assume first that there exists $\{a,b\} \in E$ such that $q \in M_{a,b} \bs \left(M_a \cup M_b\right)$. Then $M_{a,b}$ is a neighbourhood of $q$ in $M$. Since $M_{a,b}$ is locally CAT(0), $M$ is locally CAT(0) at $q$.

Assume now that there exists $a \in S$ such that $q \in M_a \bs \{q_0\}$. Then $M'_a = \bigcup_{\{a,b\} \in E} M_{a,b}$ is a neighbourhood of $q$ in $M$. Since $M'_a$ is the glueing of locally CAT(0) cube complexes $M_{a,b}$, for $\{a,b\} \in E$, along the common locally convex subspace $M_a$, we deduce that $M'_a$ is a locally CAT(0) cube complex. Therefore, $M$ is locally CAT(0) at $q$.

Assume now that $q=q_0$. We will prove that the link of $M$ at $q_0$ is flag. Assume that $Q_1,Q_2,Q_3$ are three cubes of $M$ containing $q_0$ such that each pairwise intersection $Q_i \cap Q_j$ has codimension $1$ in $Q_i$ and $Q_j$, and $Q_1 \cap Q_2 \cap Q_3$ has codimension $2$ in $Q_1$, $Q_2$ and $Q_3$. We will prove that there exists a cube $Q$ of $M$ containing $Q_1$, $Q_2$ and $Q_3$ with codimension $1$. For each $i \in \{1,2,3\}$, let $\{a_i,b_i\} \in E$ such that $Q_i \subset M_{a_i,b_i}$. 

\mk

For each $i \neq j$, we know that $M_{a_i,b_i}$ and $M_{a_j,b_j}$ intersect outside of $\{q_0\}$, so the edges $\{a_i,b_i\}$ and $\{a_j,b_j\}$ intersect in $\Gamma$. Since $\Gamma$ has no triangles, there exist $\{a,b\},\{a,c\},\{a,d\} \in E$ such that $Q_1,Q_2,Q_3 \subset M_{a,b} \cup M_{a,c} \cup M_{a,d}$. Note that $M_{a,b} \cup M_{a,c} \cup M_{a,d}$ is the union of three locally CAT(0) cube complexes along the locally convex subspace $M_a$, hence it is locally CAT(0). Therefore there exists a cube $Q$ in $M_{a,b} \cup M_{a,c} \cup M_{a,d}$ containing each of $Q_1,Q_2,Q_3$ with codimension $1$.
\ep

\section{Proofs of Corollaries} \label{sec:corollaries}

We finish by giving a proof of the two corollaries stated in the introduction

\mk

\noindent {\bf Proof of Corollary~\ref{maincor:square}}
Let $\Gamma$ be a connected, bipartite graph with diameter at least $3$ and with labels at least $3$. Since $A(\Gamma)$ has dimension $2$, the main result of~\cite{haettel_artin_cubic} applies to show that the group $A(\Gamma)$ is not virtually cocompactly cubulated. And the group $A(\Gamma)$ falls in the case (B) of Theorem~\ref{thm:main}.
\hfill \qed

\mk

\noindent {\bf Proof of Corollary~\ref{maincor:csq}}
Assume that $A$ is the fundamental group of a locally finite, finite-dimensional locally CAT(0) cube complex. Then $A$ acts freely properly by isometries on a locally finite, finite-dimensional CAT(0) cube complex.
\ben
\item According to Niblo and Reeves (see~\cite{niblo_reeves}), $A$ has the Haagerup property.
\item According to the main theorem of~\cite{guentner_higson}, $A$ is weakly amenable, with Cowling-Haagerup constant $1$.
\item According to Higson and Kasparov (see~\cite{higson_kasparov_1} and \cite{higson_kasparov_2}), the Haagerup property for $A$ implies the Baum-Connes conjecture with coefficients.
\item According to Chatterji and Ruane (see~\cite{chatterji_ruane}), since $A$ acts freely on a finite-dimensional CAT(0) cube complex, $A$ has the property RD.
\item According to Wright (see~\cite{wright}), the asymptotic dimension of a finite-dimensional CAT(0) cube complex is bounded above by its dimension.
\een
\hfill \qed

\sign

\bibliographystyle{smfalpha_perso}
\bibliography{bibli}

\end{document}